\documentclass{amsart}

\usepackage{epsfig}

\theoremstyle{plain}  
\newtheorem{thm}{Theorem}[section]
\newtheorem{prop}[thm]{Proposition}
\newtheorem{cor}[thm]{Corollary}

\theoremstyle{definition}
\newtheorem{defn}[thm]{Definition}
\newtheorem{exmp}[thm]{Example}

\theoremstyle{remark}
\newtheorem*{rem}{Remark}

\newcommand{\rank}{\mathrm{rank}}

\makeatletter 
\def\@setcopyright{} 
\def\serieslogo@{\@empty} 
\makeatother 

\begin{document}

\title[Dynkin Graphs, Gabri\'{e}lov Graphs
and Triangle Singularities]{ 
Dynkin Graphs, Gabri\'{e}lov Graphs\\
and Triangle Singularities}

\author{Tohsuke Urabe}

\address{Department of Mathematics\\
Tokyo Metropolitan University\\
Minami-Ohsawa 1-1, Hachioji-shi\\
Tokyo, 192-03 JAPAN}

\email{urabe@comp.metro-u.ac.jp}



\begin{abstract}
We consider fourteen kinds of two-dimensional triangle hypersurface 
singularities, and consider what kinds of combinations of rational double 
points can appear on small deformation fibers of these triangle singularities. 
We show that possible combinations can be described by Gabri\'{e}lov graphs
and Dynkin graphs.
\end{abstract}

\maketitle

\section{Review of results by Russian mathematicians}
\label{review}

In this article we assume that every variety is defined 
over the complex field  $\mathbf{C}$.  

First I explain some results by Arnold and Gabri\'{e}lov briefly.

In \cite{arnold;gauss} Arnold has introduced an invariant  $m$  called 
\emph{modality} or \emph{modules 
number}, and has given a long classification list of hypersurface 
singularities.  
Modality $m$ is a non-negative integer.
Though we find singularities of any dimension in Arnold's list, 
we consider singularities of dimension two in particular.  

His class of singularities with  $m=0$ coincides with the class of rational 
double points.  
It is well known that each rational double point corresponds to a 
connected Dynkin graph of type  $A$, $D$  or  $E$  in the theory of Lie 
algebras.  (Durfee~\cite{durfee;charac}.)  

The class with $m=1$ consists of three subclasses.  
($\lambda$  is a parameter.)

\begin{enumerate}
\item Three simple elliptic singularities: $J_{10},\;X_9,\;P_8$
\item Cusp singularities $T_{p,\,q,\,r}$.
$\left( {\frac{1}{p}+\frac{1}{q}+\frac{1}{r}<1} \right)$:
$x^p+y^q+z^r+\lambda xyz=0$ $\left( {\lambda \ne 0} \right)$.
\item fourteen triangle singularities 
(These fourteen are also called exceptional singularities.)
$$\begin{array}{llllllll}
E_{12}&Z_{11}&Q_{10}&\hspace{10 mm}&W_{12}&S_{11}&\hspace{10 mm}&U_{\,12}\\
E_{13}&Z_{12}&Q_{11}&\hspace{10 mm}&W_{13}&S_{12}\\
E_{14}&Z_{13}&Q_{12}
\end{array}$$
$$\begin{array}{c}
E_{12}\,:\;x^7+y^3+z^2+\lambda x^5y=0\hspace{5 mm}
W_{12}\,:\;x^5+y^4+z^2+\lambda x^3y^2=0\\
U_{12}\,:\;x^4+y^3+z^3+\lambda x^2yz=0.
\end{array}$$
\end{enumerate}
(As for the other defining polynomials see Arnold~\cite{arnold;gauss}.)

His list continues in the case $m\ge 2$,
but we do not refer further.  

We go on to Gabri\'{e}lov's results.  
(Gabri\'{e}lov~\cite{gabrielov;unimodular}.)

Let $f\left( {x,\,y,\,z} \right)=0$
be one of defining polynomials of fourteen hypersurface triangle 
singularities.  
It defines a singularity at the origin.  
We consider the \emph{Milnor fiber}, i.e.,
$$F=\left\{ {\,\left( {x,\,y,\,z} \right)\in \mathbf{C}^3\;\left| 
{\;\left| x \right|^2+\left| y \right|^2+\left| z \right|^2<\varepsilon ^2,
\;f\left( {x,\,y,\,z} \right)=t\,} \right.} \right\}$$
where  $\epsilon$  is a sufficiently small positive real number and  $t$  is 
a non-zero complex number whose absolute value is sufficiently small 
compared with  $\epsilon$.  
The pair
$$
\left( H_2\left( F,\;\mathbf{Z} \right),\;the\ intersection\ form\right)
$$
is called the \emph{Milnor lattice}, 
and $\mu =\rank \,H_2\left( F,\;\mathbf{Z} \right)$ is called the 
\emph{Milnor number} of the singularity.  
Gabri\'{e}lov has computed the Milnor lattice for fourteen hypersurface 
triangle singularities.  
According to him, there exists a basis $e_{\,1},\;e_2,\;\ldots ,\;e_\mu $
of $H_2\left( F,\;\mathbf{Z} \right)$
such that each $e_{\,i}$ is a vanishing cycle 
(thus in particular $e_{\,i}\cdot e_{\,i}=-2$)  
and the intersection form is represented by the dual graph below. 

\vspace*{2 mm}
\epsfig{file=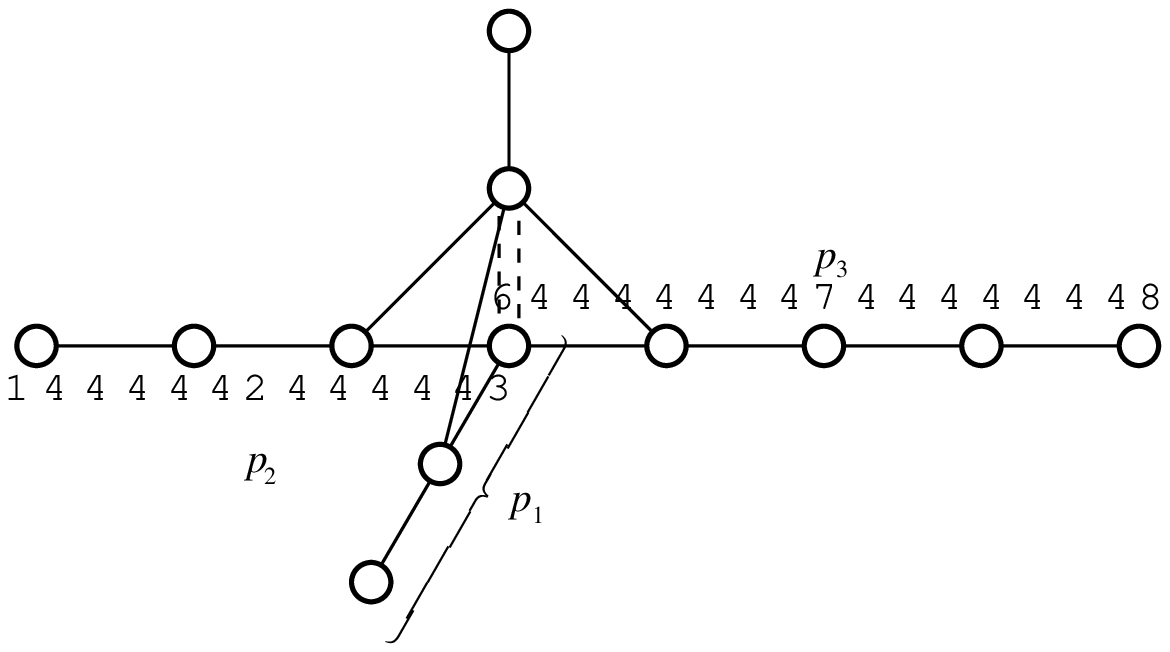}

In the above graph the basis $e_{\,1},\;e_2,\;\ldots ,\;e_\mu $
has one-to-one correspondence with vertices.  
Edges indicate intersection numbers.  
Two vertices corresponding to  $e_{\,i}$
and  $e_j$ are not connected, if $e_{\,i}\cdot e_j=0$.
They are connected by a single solid edge, if $e_{\,i}\cdot e_j=1$.
They are connected by a double dotted edge, if $e_{\,i}\cdot e_j=-2$. 
Three integers $p_1,\;p_2,\;p_3$ are the numbers of vertices in the 
corresponding arm.\nopagebreak[4]
They depend on the type of the triangle singularity. 
The corresponding triplets $\left( {p_{\,1},\;p_2,\;p_3} \right)$
to the above fourteen symbols are as follows:
$$\begin{array}{llllllll}
\left( {2,\;3,\;7} \right)&\left( {2,\;4,\;5} \right)&
\left( {3,\;3,\;4} \right)&\hspace{10 mm}&
\left( {2,\;5,\;5} \right)&\left( {3,\;4,\;4} \right)&\hspace{10 mm}&
\left( {4,\;4,\;4} \right)\\
\left( {2,\;3,\;8} \right)&\left( {2,\;4,\;6} \right)&
\left( {3,\;3,\;5} \right)&\hspace{10 mm}&
\left( {2,\;5,\;6} \right)&\left( {3,\;4,\;5} \right)\\
\left( {2,\;3,\;9} \right)&\left( {2,\;4,\;7} \right)&
\left( {3,\;3,\;6} \right)
\end{array}$$

(Thus the above figure is the graph for $S_{12}$.)  

The main part of the above graph below is called the 
\emph{Gabri\'{e}lov graph}.  

\vspace*{2 mm}
\epsfig{file=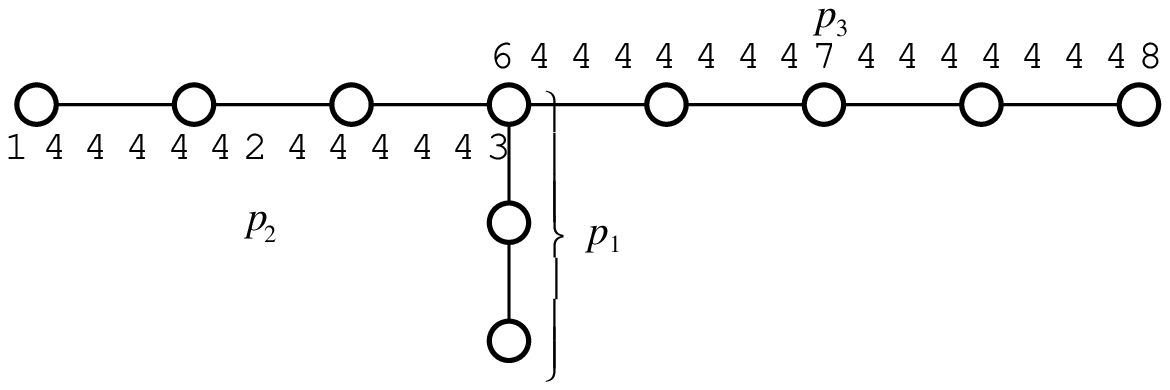}

\vspace*{\fill}
\pagebreak[4]
The Gabri\'{e}lov graph defines a lattice  $P^*$  
with a basis $e_{\,1},\;e_2,\;\ldots ,\;e_{\mu -2}$
, if we apply the above mentioned rule.  
It is easy to check that  $P^*$  has signature $\left( 1,\;\mu -3 \right)$,
and $H_2\left( F,\;\mathbf{Z} \right)\cong P^*\oplus H$ as lattices.  
Here $H=\mathbf{Z}\,u+\mathbf{Z}\,v$ denotes the hyperbolic plane, i.e., 
a lattice of rank 2 
with a basis  $u$, $v$  satisfying $u\cdot u=v\cdot v=0$ and 
$u\cdot v=v\cdot u=1$,
and $\oplus $ denotes the orthogonal direct sum.  

\nopagebreak[4]
\section{Singularities on deformation fibers of triangle singularities}
\label{triangle}

A finite disjoint union of connected Dynkin graphs is also called a Dynkin 
graph.  

Let  $T$  denote one of the above fourteen symbols of hypersurface triangle 
singularities.  
By $PC\left( T \right)$ we denote the set of Dynkin graphs  $G$  with 
several components such that there exists a small deformation fiber  $Y$  
of a singularity of type  $T$  satisfying the following conditions:
\begin{enumerate}
\item  $Y$  has only rational double points as singularities.
\item The combination of rational double points on  $Y$  corresponds to 
graph  $G$  exactly.  
\end{enumerate}
Here, the type of each component of  $G$  corresponds to the type of a 
rational double point on  $Y$,  and the number of components of each type 
corresponds to the number of rational double points of each type on  $Y$.
If $G$  has $a_k$ of components of type $A_k$ for each $k\ge 1$,
$d_\ell $ of components of type $D_\ell $ for each  $\ell \ge 4$ and  
$e_m$  of components of type  $E_m$  for  $m=6,\;7,\;8$,
we identify  $G$  with the formal sum  $G=\sum {a_k\,A_k}+\sum {d_\ell \,
D_\ell} +\sum {e_m\,E_m}$.

Mr. F.-J. Bilitewski informed me that he had a complete listing of Dynkin 
graphs of  $PC\left( T \right)$ for every  $T$  of the above fourteen.  

\begin{thm}
\label{main}
 Let  $T$  be one of the above fourteen symbols of hypersurface triangle 
singularities.  
Let  $G$  be a Dynkin graph with only components of type  $A$, $D$  or  $E$.  
The following conditions \textup{\textbf{(A)}} and \textup{\textbf{(B)}} are 
equivalent:
\begin{description}
\item[(A)] $G\in PC\left( T \right)$.
\item[(B)] Either \textup{\textbf{(B-1)}} or \textup{\textbf{(B-2)}} holds.
\begin{description}
\item[(B-1)]   $G$  can be made by an elementary transformation or a tie 
transformation from a Dynkin subgraph of the Gabri\'{e}lov graph of type  
$T$.  
\item[(B-2)]  $G$  is one of the following exceptions:
\end{description}
\end{description}

\vspace*{-5 mm}
\[
\begin{array}{l}
\begin{array}{rl}
 &\mbox{\hspace{5 mm}Exceptions}\\
T & =Z_{\,13}\,:\;A_7+A_4\\
T & =S_{\,11}\,:\;2A_4+A_{\,1}\\
T & =U_{\,12}\,:\;2D_4+A_2,\;A_6+A_4,\;A_5+A_4+A_{\,1},\;2A_4+A_1
\end{array}
\\
\mbox{The other eleven triangle singularities: None}
\end{array}
\]
\end{thm}

An elementary transformation and a tie transformation in the above 
theorem are operations by which we can make a new Dynkin graph from a 
given Dynkin graph.  

\begin{defn}
Elementary transformation:  
The following procedure is called an elementary transformation of a 
Dynkin graph:
\begin{enumerate}
\item Replace each connected component by the corresponding 
extended Dynkin graph.
\item Choose in an arbitrary manner at least one vertex from each 
component (of the extended Dynkin graph) and then remove these vertices 
together with the edges issuing from them.
\end{enumerate}
\end{defn}

An extended Dynkin graph is a graph obtained from a connected Dynkin graph 
by adding one vertex and one or two edges.  (Bourbaki~\cite{bourbaki;lie}.)  
Below we show extended Dynkin graphs.  
Numbers in the figures below are the coefficients of the maximal root, 
which will appear in the definition of a tie transformation.  
We can get the corresponding Dynkin graph, 
if we erase one vertex with the attached number $1$ and edges issuing 
from it.  

\vspace*{5 mm}
\epsfig{file=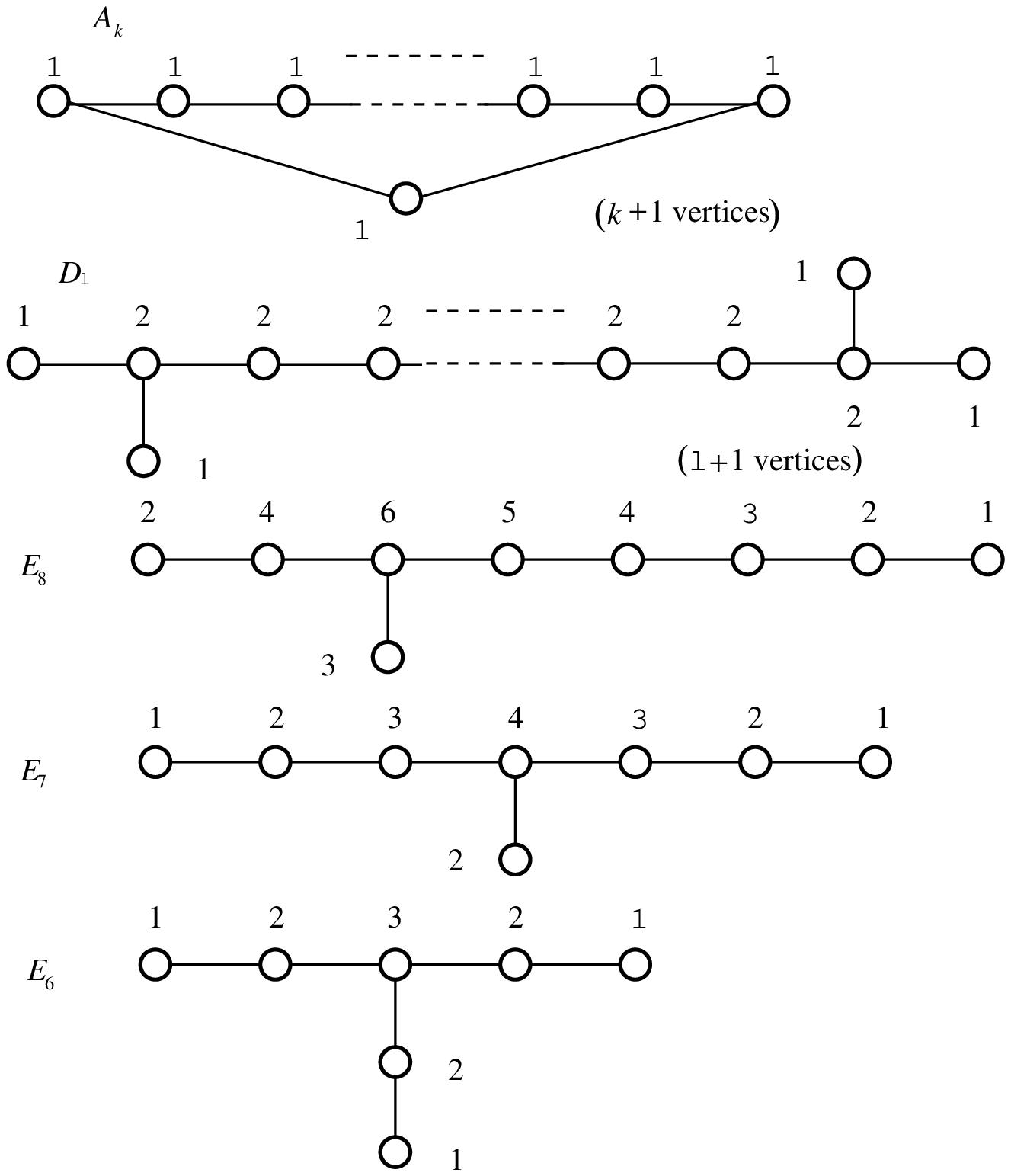}

\vspace{\fill}

\pagebreak[4]

\begin{defn}
 Tie transformation:  
Assume that by applying the following procedure to a Dynkin graph  $G$  
we have obtained the Dynkin graph $\bar G$.
Then, we call the following procedure a tie transformation of a 
Dynkin graph:
\begin{enumerate}
\item Replace each component of  $G$  by the extended Dynkin graph of the 
same type.  
Attach  the corresponding coefficient of the maximal root to each vertex of 
the resulting extended graph $\tilde G$.
\item Choose, in an arbitrary manner, subsets  $A$, $B$  of the set of 
vertices of the extended graph $\tilde G$ satisfying the following 
conditions:
\begin{description}
\item[$\left\langle a \right\rangle$] $A\cap B=\emptyset $
\item[$\left\langle b \right\rangle$] 
Choose arbitrarily a component  $\tilde G''$ of $\tilde G$ and let  $V$  be 
the set of vertices in  $\tilde G''$.  
Let  $\ell $ be the number of elements in  $A\cap V$.  
Let $n_{\,1},\;n_2,\;\ldots ,\;n_\ell $ be the numbers attached to 
$A\cap V$.
Also, let  $N$  be the sum of the numbers attached to elements in 
$B\cap V$. (If  $B\cap V=\emptyset $, $N=0$.)  
Then, the greatest common divisor of the  $\ell +1$ numbers 
$N,\;n_{\,1},\;n_2,\;\ldots ,\;n_\ell$ is $1$.
\end{description}
\item Erase all attached integers.
\item Remove vertices belonging to  $A$  together with the edges issuing 
from them.
\item Draw another new vertex called $\theta $.
Connect  $\theta $ and each vertex in  $B$  by a single edge.
\end{enumerate}
\end{defn}

\begin{rem}
After following the above procedure 1--5, the resulting graph  $\bar G$
is often not a Dynkin graph.  
We consider only the cases where the resulting graph  $\bar G$  is a Dynkin 
graph, and then we call the above procedure a tie transformation.
The number  $\#\left( B \right)$ of elements in the set  $B$  satisfies  
$0\le \#\left( B \right)\le 3$.
$\ell =\#\left( {A\cap V} \right)\ge 1$.
\end{rem}

\begin{exmp}
We consider the case $T=W_{\,13}$.
The Gabri\'{e}lov graph in this case is the following, and 
it has a Dynkin subgraph of type $E_8+A_2$:

\vspace*{3 mm}
\epsfig{file=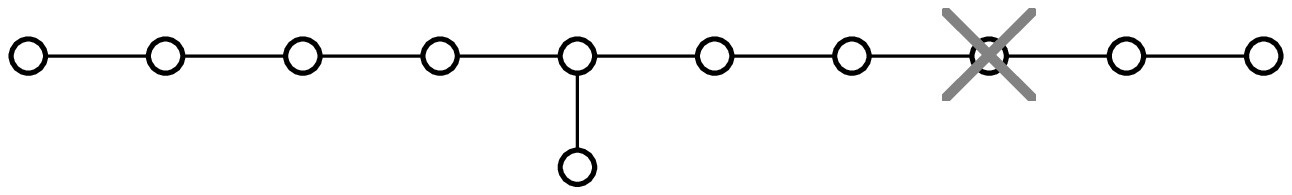}
\vspace*{3 mm}

First we apply a tie transformation to  $E_8+A_2$.
In the second step of the transformation we can choose subsets  $A$  
and  $B$  as follows:

\vspace*{4 mm}
\epsfig{file=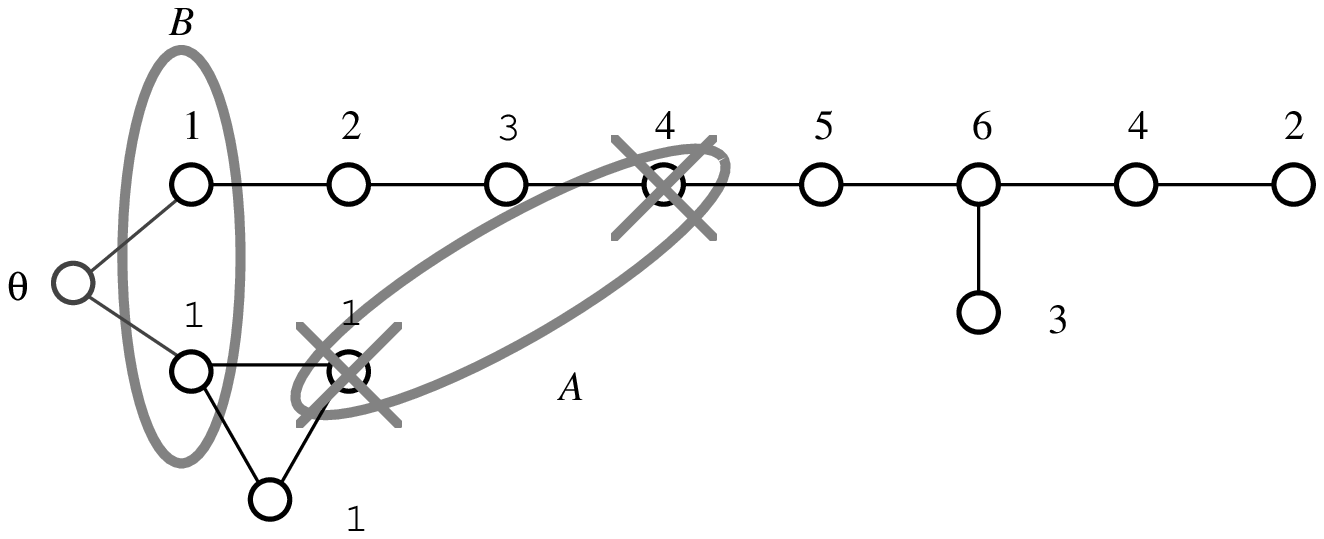}

\pagebreak[4]
For the component of type $E_8$, $\ell =1$, $n_{\,1}=4$, $N=1$
and thus $G.C.D.\left( {n_{\,1},\;N} \right)=1$.
For the component $A_{\,2}$, $\ell =1$, $n_{\,1}=1$, $N=1$ and 
thus $G.C.D.\left( {n_{\,1},\;N} \right)=1$. 
One sees that the condition $\left\langle b \right\rangle $ is satisfied.  
As the result of the transformation one gets a graph of type $A_6+D_5$.
By our theorem one can conclude $A_6+D_5\in PC\left( {W_{\,13}} \right)$.

Second we apply an elementary transformation to $E_8+A_2$.

\vspace*{4 mm}
\epsfig{file=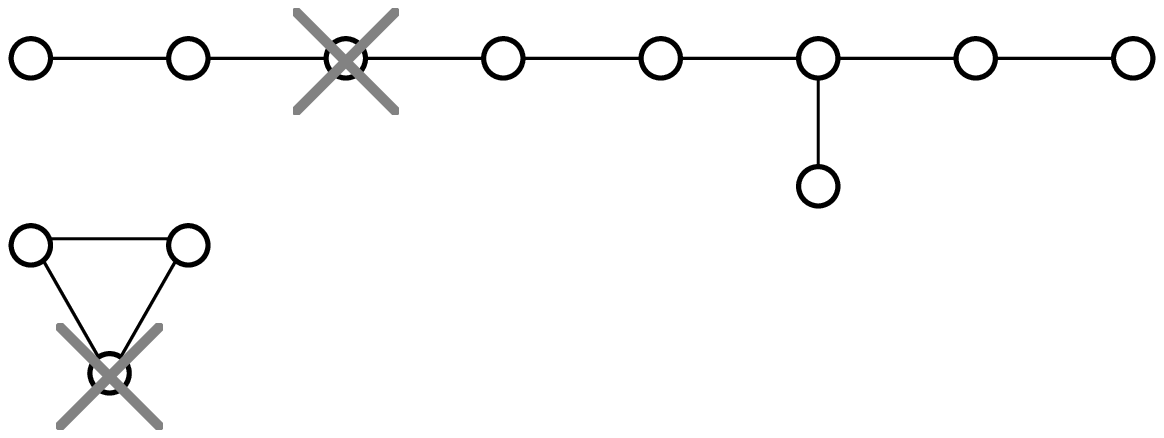}
\vspace*{3 mm}

As in the above figure we can get $E_6+2A_2$.
Thus $E_6+2A_2\in PC\left( {W_{\,13}} \right)$.
\end{exmp}

\nopagebreak[4]
\section{K3 surfaces and lattice theory}
\label{k3}

It is known that fourteen hypersurface triangle singularities have 
interesting property called the strange duality.  
(Pinkham~\cite{pinkham;duality}.)  
Let  $T$  be one of the above fourteen symbols of hypersurface triangle 
singularities.  
Associated with  $T$, we have another symbol  $T^*$ also in the above 
fourteen symbols of hypersurface triangle singularities.  
This  $T^*$  is called the \emph{dual} of  $T$.  
The dual of the dual coincides with the original one, i.e.,  $(T^*)^* = T$.  

For the following six singularities the dual coincides with itself, i.e.,  
$T^* = T$: $E_{\,12},\;Z_{\,12},\;Q_{\,12},\;W_{\,12},\;S_{\,12},\;U_{\,12}$.
For the following four pairs the dual is another member of the pair: 
$\left\{ {\,E_{\,13},\;Z_{\,11}\,} \right\},\;\left\{ {\,E_{\,14},\;
Q_{\,10}\,} \right\},\;\left\{ {\,Z_{\,13},\;Q_{\,11}\,} \right\},\;
\left\{ {\,W_{\,13},\;S_{\,11}\,} \right\}$.

Following Looijenga~\cite{looijenga;triangle}, 
we explain the relation between triangle singularities and K3 surfaces 
below.  
Let  $T$  be one of the above fourteen symbols of hypersurface triangle 
singularities.  
Let  $\Gamma^*$ be  the Gabri\'{e}lov graph of the dual  $T^*$.  
We can define a reducible curve  $IF$  on a surface whose dual graph 
coincides with $\Gamma^*$.  
The curve $IF=IF\left( T \right)$ is called \emph{the curve at infinity} of type  $T$. 
The irreducible components are all  smooth rational curves  $C$  
with $C\cdot C=-2$ and have one-to-one correspondence with vertices of  
$\Gamma^*$.  
For two components  $C$, $C'$ of  $IF$  the intersection number $C\cdot C'$
is equal to one or zero, according as the corresponding vertices in  $\Gamma^*$  
are connected in  $\Gamma^*$  or not.  

Let  $G$  be a Dynkin graph with components of type  $A$, $D$  or  $E$  only.  
Assume that there exists a smooth K3 surface  $Z$  satisfying the 
following conditions \textbf{(a)} and \textbf{(b)}:
\begin{description}
\item[(a)] $Z$  contains the curve at infinity $IF=IF\left( T \right)$ of type  
$T$  as a subvariety.  
\item[(b)] Let  $E$  be the union of all smooth rational curves on  $Z$  
disjoint from  $IF$.  
The dual graph of the components of  $E$  coincides with graph  $G$.  
\end{description}
(Note that an irreducible curve  $C$  on a K3 surface is a smooth rational 
curve if, and only if, $C\cdot C=-2$.)  

Contracting every connected component of  $E$  to a rational double point 
and then removing the image of  $IF$, we obtain an open variety $\tilde Y$.

\begin{prop}[Looijenga~\cite{looijenga;triangle}]
\label{looijenga}
\begin{enumerate}
\item  Under the above assumption there exists a small deformation fiber  
$Y$  of a singularity of type  $T$  homeomorphic to $\tilde Y$.  
\item Let  $Y$  be a small deformation fiber of a singularity of type  $T$.  
Assume that  $Y$  has only rational double points as singularities, 
and the combination of rational double points on  $Y$  corresponds to a 
Dynkin graph  $G$.  
Then, there exists a K3 surface satisfying \textup{\textbf{(a)}} and 
\textup{\textbf{(b)}}, 
and the corresponding $\tilde Y$ is homeomorphic to  $Y$.  
\end{enumerate}
\end{prop}

By the above proposition our study is reduced to the study of K3 surfaces 
containing the curve $IF=IF\left( T \right)$.
K3 surfaces are complicated objects, but it is known that by the theory of 
periods we can reduce the study of K3 surfaces to the study of lattices.  

Below we explain several terminologies in the lattice theory.  
(Urabe~\cite{urabe;quardli}.)  
A free module over  $\mathbf{Z}$  of finite rank equipped with an integral 
symmetric bilinear form $\left( {\;\;\;,\;\;} \right)$ is called a 
\emph{lattice}.  
Besides, if a free module  $L$  over  $\mathbf{Z}$  of finite rank has a 
symmetric bilinear form $\left( {\;\;\;,\;\;} \right)$ with values in rational 
numbers, then  $L$  is called a \emph{quasi-lattice}.  
For simplicity we write $x^2=\left( {x,\;x} \right)$.

Let  $L$  be a quasi-lattice and  $M$  be a submodule.  
The submodule 
$$\tilde M=\left\{ x\in L\, \left| {\,mx\in M\ } 
\mbox{for some non-zero integer} \ m\,\right. \right\}$$
 is called the \emph{primitive hull} of  $M$  in  $L$.  
We say that  $M$  is \emph{primitive}, if $M=\tilde M$, 
and an element $x\in L$ is primitive, 
if $M=\mathbf{Z}\,x$ is primitive.  
We say that an embedding $M\to L$ of quasi-lattices is a \emph{primitive 
embedding}, if the image is primitive.  
If  $M$  is  non-degenerate and primitive as a sub-quasi-lattice, 
we can define the canonical induced bilinear form on the quotient module
$L/M$.

Let  $L$  be a quasi-lattice, and  $FL$  be a submodule of  $L$  such that 
the index $\#\left( L/ FL \right)$ is finite.  Set
\begin{eqnarray*}
R=\left\{ \,\alpha \in FL\;\left| \;\alpha ^2=-2\, \right. \right\}
\cup\left\{ \,\beta \in L\;\left|
\;\beta ^2=-1\ or\ -2/ 3\; 
\right.\right\}\\
\hspace{7 mm}\cup  \left\{ \,\gamma \in L\;\left| \;\gamma ^2=-1/ 2,\;
2\gamma \in FL\; \right. \right\}
\end{eqnarray*}

The set $R=R(L,\;FL)$ is called the \emph{root system} of 
$\left( {L,\;FL} \right)$, 
and every element $\alpha \in R$ is called a \emph{root}.  
If the pair $\left( {L,\;FL} \right)$ satisfies the following conditions 
\textbf{(R1)} and \textbf{(R2)}, then $\left( {L,\;FL} \right)$ is called a 
\emph{root module}:
\begin{description}
\item[(R1)] $2\left( {x,\;\alpha } \right)/\alpha ^2$ is an integer for 
every $x\in L$ and $\alpha \in R$.
\end{description}

Under \textbf{(R1)}, for every $\alpha \in R$ we can define an isomorphism  
$s_\alpha \,:\;L\to L$ preserving the bilinear form, by setting for $x\in L$ 
$s_\alpha \left( x \right)=x-2\left( {x,\;\alpha } \right)\alpha / 
\alpha ^2$.
\begin{description}
\item[(R2)]$s_\alpha \left( {FL} \right)=FL$ for every $\alpha \in R$.
\end{description}

Let $\left( {L,\;FL} \right)$ be a root module. 
If $L=FL$,  we say that it is \emph{regular} and abbreviate  $FL$.   
Let  $M$  be a submodule of  $L$.  
It is easy to check that the pair $\left( {M,\;FL\cap M} \right)$ is again a 
root module.  
Below we identify  $M$  with the pair $\left( {M,\;FL\cap M} \right)$.  
If the root system of  $M$  and the root system of $\tilde M$ coincide, 
then we say that  $M$  is \emph{full}.  
An embedding $M\to L$ of quasi-lattices is a \emph{full embedding}, 
if the image is full.  

Let  $G$  be a Dynkin graph with several components of type $A$, $D$ or  
$E$ only.  
We can define a lattice and its basis such that the corresponding dual 
graph coincides with  $G$.  
This lattice is called the \emph{root lattic}e of type  $G$  and is denoted 
by $Q(G)$. 
$Q(G)$ is a regular root module with a basis $\alpha _{\,1},\;\alpha _{\,2},\;
\ldots ,\;\alpha _r$ with $\alpha _{\,i}^2=-2$ for every  $i$.  
Let $\Lambda _N$ denote the even unimodular lattice with 
signature $\left( {N,\;16+N} \right)$ for $N\ge 0$.
The isomorphism class of $\Lambda _N$ is unique if $N\ge 1$ and 
thus $\Lambda _N\cong \Lambda _{N-1}\oplus H$.
For a K3 surface  $Z$  the second cohomology group 
$H^2\left( {Z,\;\mathbf{Z}} \right)$ with the intersection form is a lattice 
isomorphic to $\Lambda _3$. 
Let $P=P\left( T \right)$ be the lattice whose dual graph is the Gabri\'{e}lov 
graph  $\Gamma^*$  of the dual  $T^*$.  
Assume that there exists a K3 surface  $Z$  satisfying the above condition 
\textbf{(a)}.  
The classes of the components of $IF$ generate a primitive sublattice in 
$H^2\left( {Z,\;\mathbf{Z}} \right)$,  which is isomorphic to  $P$.

\begin{prop}
\label{embedding}
\begin{enumerate}
\item  If $N\ge 1$, there is a primitive embedding $P\to \Lambda _N$.
\item  If $N\ge 2$,  a primitive embedding $P\to \Lambda _N$ is unique up 
to automorphisms of $\Lambda _N$.
\item If $N\ge 1$, for any embedding $P\to \Lambda _N$,  
the pair $\left( {\Lambda _N/\tilde P,\;F_N} \right)$ is a root module, 
where  $F_N$ is the image of the orthogonal complement of  $P$  
in $\Lambda _N$ by the canonical surjective homomorphism
$\Lambda _N\to \Lambda _N/ \tilde P$.
\item  For any primitive embedding $P=P\left( T \right)\to \Lambda _2$ the 
orthogonal complement $F_2$ of  $P$  in  $\Lambda _2$ has a basis whose 
dual graph coincides with the Gabri\'{e}lov graph of type $T$.  
\end{enumerate}
\end{prop}

With aid of Looijenga's results in \cite{looijenga;triangle} 
we can show the following:

\begin{prop}
\label{full}
We fix a primitive embedding $P\to \Lambda _3$.  
There exists a K3 surface  $Z$  satisfying the above conditions 
\textup{\textbf{(a)}} and \textup{\textbf{(b)}} if, and only if, 
there is a full embedding $Q\left( G \right)\to \Lambda _3/P$.
\end{prop}

\begin{cor}
\label{criterion}
 $G\in PC\left( T \right)$  if, and only if, there is a full embedding  
$Q\left( G \right)\to \Lambda _3/P\left( T \right)$.
\end{cor}

By Proposition~\ref{full} our study has been reduced to the lattice theory.  
Next, we have to consider properties of the lattice $P=P\left( T \right)$ 
depending on  $T$  closely.  
Let  $T$  be one of fourteen symbols of hypersurface triangle singularities.  

\begin{prop}
\label{property}
 We fix $N\ge 1$.
\begin{enumerate}
\item  For any  $T$  and for any embedding  $P(t)\to \Lambda _N$ the 
quasi-lattice  $\Lambda _N/\tilde P\left( T \right)$ does not contain an 
element $\beta$  with  $\beta ^2=-1$.
\item The root module $\left( {\Lambda _N/\tilde P\left( T \right),\;F_N} 
\right)$ contains a root $\gamma$ with $\gamma ^2=-1/2$ for some 
embedding $P\left( T \right)\to \Lambda _N$ if, and only if, 
$T=E_{\,13},\;Z_{\,12},\;Q_{\,11},\;W_{\,13}$ or $U_{\,12}$. 
It contains a root $\gamma$ with $\gamma ^2=-1/2$ for some primitive 
embedding $P\left( T \right)\to \Lambda _N$ if, and only if, 
$T=E_{\,13},Z_{\,12}$ or $\ Q_{\,11}$.
\item The root module $\left( {\Lambda _N/\tilde P\left( T \right),\;F_N} 
\right)$ contains a root $\beta$ with $\beta ^2=-2/3$ for some embedding 
$P\left( T \right)\to \Lambda _N$ if, and only if, $T=E_{\,14},\;Z_{\,13}$
or $Q_{\,12}$.
\end{enumerate}
\end{prop}

Consider the case where $(L,\;FL)$ is a root module such that the bilinear
form on $L$ has signature $(1,\;\rank\, L-1)$.
In this case we can apply the hyperbolic geometry, and we can give the
generalization of the theory in the negative definite case such as the Weyl
chamber and the Dynkin graph.
The generalized Dynkin graph in this case is called the 
\emph{Coxeter-Vinberg graph}. 
(Vinberg~\cite{vinberg;group}.)

We need consider the Coxeter-Vinberg graph of 
$\left( {\Lambda _{\,2}/P\left( T \right),\;F_2} \right)$.  
By Proposition~\ref{embedding}.4 we can expect that it is related to the 
Gabri\'{e}lov graph.  
We fix a primitive embedding $P\left( T \right)\to \Lambda _2$.

\begin{prop}
\label{vinberg-g}
 Let $\tilde \Gamma$ denote the Coxeter-Vinberg graph of 
$\left( {\Lambda _{\,2}/P\left( T \right),\;F_2} \right)$.
\begin{enumerate}
\item  We can draw $\tilde \Gamma$ in finite steps if, and only if, 
$T\ne S_{\,11},\;S_{\,12}$.
\item  If $T\ne W_{\,12},\;W_{\,13},\;S_{\,11},\;S_{\,12},\;\ U_{\,12}$,  
every vertex in $\tilde \Gamma$ corresponds to a root.  
\item  If $T=W_{\,12}$,  every vertex in $\tilde \Gamma$ corresponds 
to either a root $\alpha$ with $\alpha ^2=-2$ or an element $\delta$ 
with $\delta ^2=-2/ 5$.
\item  If $T=W_{\,13}\ $ or $U_{\,12}$,  every vertex in $\tilde \Gamma$ 
corresponds to either a root $\alpha$  with $\alpha ^2=-2$ 
or an element $\delta$ with $\delta ^2=-1/2$ and $2\delta \notin F_2$.
\item If $T=E_{\,12},\;Z_{\,11},$ or $\ Q_{\,10}$,  the Gabri\'{e}lov graph 
coincides with $\tilde \Gamma$.
\item  If $T=E_{\,13},\;E_{\,14},\;Z_{\,12},\;Z_{\,13},\;Q_{\,11},\;Q_{\,12},
\;W_{\,12},$ or $\ U_{\,12}$,  the Gabri\'{e}lov graph is the subgraph of 
$\tilde \Gamma$ consisting of all vertices corresponding to a root $\alpha$
with $\alpha ^2=-2$.
\item  If $T=W_{\,13}$,  the Gabri\'{e}lov graph is the maximal subgraph of  
$\tilde \Gamma$ such that every vertex corresponds to a root $\alpha$ with 
$\alpha ^2=-2$,  and  if $\alpha$, $\beta$ are roots corresponding to two 
vertices, then $\left( {\alpha ,\;\beta } \right)\ne -2$.
\end{enumerate}
\end{prop}

We can explain main ideas in the verification of our Theorem~\ref{main} 
here.  
Let $\overline {PC}\left( T \right)$ denote the set of all Dynkin graphs 
made from a Dynkin subgraph of the Gabri\'{e}lov graph of type  $T$  
by an elementary transformation or a tie transformation.  
We assume that $G\in \overline {PC}\left( T \right)$ was made from a 
Dynkin subgraph $G'$ of the Gabri\'{e}lov graph.  
Besides, we fix a primitive embedding $P\to \Lambda _N$ 
for $N=2,\;3$.  
By Proposition~\ref{embedding}.4 there is a primitive embedding  
$Q\left( {G'} \right)\to F_2$. 
By the theory of elementary and tie transformations 
(Urabe~\cite{urabe;elem}, \cite{urabe;tie}.) we can conclude that there is 
a full embedding $Q\left( G \right)\to F_3\cong F_2\oplus H$ into the 
regular root module $F_3$.  
Assume $T\ne E_{\,13},\;Z_{\,12},\;Q_{\,11},\;E_{\,14},\;Z_{\,13},
\;Q_{\,12}$ here.  
By Proposition~\ref{property} the composition 
$Q\left( G \right)\to F_3\subset \Lambda _3/P$ defines a full embedding 
into the root module $\left( {\Lambda _3/P,\;F_3} \right)$ in these cases.  
By Corollary~\ref{criterion} we have $G\in PC\left( T \right)$.  
Thus $\overline {PC}\left( T \right)\subset PC\left( T \right)$.

Next, we determine the difference 
$PC\left( T \right)-\overline {PC}\left( T \right)$.  
Let  $r$  be the number of vertices of a graph  $G$.  
It is easy to see that if $G\in PC\left( T \right)$,  then $r\le \mu -2$.  
In case $T\ne S_{\,11},\;S_{\,12}$,  using Proposition~\ref{vinberg-g} 
we can show that conditions $G\in PC\left( T \right)$ and  
$G\in \overline {PC}\left( T \right)$ are equivalent if $r\le \mu -5$.  
Thus we can assume $r=\mu -2,\;\mu -3$ or $\mu -4$.  
For triangle singularities the Milnor number $\mu$ is relatively small, 
and it is easy to check whether a Dynkin graph  $G$  belongs to  
$PC\left( T \right)-\overline {PC}\left( T \right)$ case-by-case.  
To tell the truth, we could not succeed in finding any effective method 
except case-by-case checking.  
This is a weak point of our theory.  
I regret this fact and hope that somebody can improve it.  
If $T=S_{\,11}$ or $S_{\,12}$, the checking becomes more complicated 
since we have no Coxeter-Vinberg graph.  

Now, if $T\ne W_{\,12},\;W_{\,13},\;S_{\,11},\;S_{\,12},\;\ U_{\,12}$,  
then because of Proposition~\ref{vinberg-g}.2 we can formulate another 
theorem.  
(Urabe~\cite{urabe;triangle}.) 
In this another theorem we start from not a Gabri\'{e}lov graph but a 
Dynkin graph possibly with a component of type $BC_{\,1}$ or $G_2$, 
and the number of transformations is not one but two.  
There, no exception appears even in the case $Z_{\,13}$.  
(We can make $A_7+A_4$ from $E_7+G_2$ by two tie transformations.)  
For $T=E_{\,13},\;Z_{\,12},\;Q_{\,11},\;E_{\,14},\;Z_{\,13}$ or 
$Q_{\,12}$ our theorem in this article follows from this theorem 
in another formulation.  
Also for $T=W_{\,12},\;W_{\,13},\;S_{\,11},\;S_{\,12},\;U_{\,12}$ 
the theorem in another formulation is possible, but becomes very 
complicated, because Proposition~\ref{vinberg-g}.2 does not hold for them.  
It is not worth mentioning.  

Details of the verification will appear elsewhere.

Now, it is very strange that our Theorem~\ref{main} has a few exceptions 
in a few cases.  
Perhaps this is because our theory has a missing part. 

\vspace*{3 mm}
\noindent
\textbf{Problem.}  Find the missing part of our theory and give a simple 
characterization of the set $PC\left( T \right)$ without exceptions.  

\vspace*{3 mm}
This problem may be very difficult, but I believe that there exists a solution.

\end{document}